\documentclass[a4paper,10pt]{scrartcl}
\usepackage{amsmath,amsthm}
\usepackage{amsfonts}
\usepackage{amssymb}
\usepackage{a4wide}
\usepackage{amscd}
\usepackage{color}
\usepackage{graphicx}
\usepackage{float}
\usepackage[utf8]{inputenc}
\usepackage{cleveref}
\usepackage{algorithm}
\usepackage{algpseudocode}

\usepackage{multirow}

\newcommand{\N}{\mathbb{N}}
\newcommand{\R}{\mathbb{R}}

\newcommand{\xdag}{x^\dag}

\newcommand{\supp}{\mathrm{supp}}
\newcommand{\csupp}{\mathrm{conv\,supp}}

\newenvironment{frcseries}{\fontfamily{frc}\selectfont}{}

\DeclareMathOperator*{\argmin}{arg\,min}

\numberwithin{equation}{section}

\newtheoremstyle{slanted}
{}
{}
{\slshape}
{}
{\bfseries}
{.}
{ }
{}
\newtheoremstyle{roman}{}{}{\rmfamily}{}{\bfseries}{.}{ }{}

\theoremstyle{plain}
\newtheorem{theorem}{Theorem}
\newtheorem{definition}{Definition}
\newtheorem{lemma}{Lemma}

\theoremstyle{slanted}

\theoremstyle{roman}

\title{On uniqueness and ill-posedness for the deautoconvolution problem in the multi-dimensional case}
\author{\textsc{Bernd Hofmann}\thanks{Faculty of Mathematics, Chemnitz University of Technology, 09107 Chemnitz, Germany, \newline e-mail: hofmannb@mathematik.tu-chemnitz.de.},\quad\textsc{Frank Werner}\thanks{Institute for Mathematics, University of W\"urzburg, Emil-Fischer-Str.~30, 97074~W\"urzburg, Germany, e-mail: frank.werner@mathematik.uni-wuerzburg.de} \,\,and
\textsc{Yu Deng}\thanks{Chemnitz University of Technology, Faculty of Mathematics, 09107 Chemnitz, Germany, e-mail: yu.deng@math.tu-chemnitz.de}}

\date{\today}

\begin{document}

\maketitle

\begin{quote}
\textbf{Abstract:} This paper analyzes the inverse problem of deautoconvolution in the
multi-dimensional case with respect to solution uniqueness and ill-posedness.
Deautoconvolution means here the reconstruction of a real-valued $L^2$-function
with support in the $n$-dimensional unit cube $[0,1]^n$ from observations of its
autoconvolution either in the full data case (i.e. on $[0,2]^n$) or in the limited data case (i.e. on $[0,1]^n$). Based on multi-dimensional variants of the Titchmarsh
convolution theorem due to Lions and Mikusi\'{n}ski, we prove in the full data case a twofoldness assertion, and in the limited data case uniqueness of non-negative solutions for which the origin belongs to the support. The latter assumption is also shown to be necessary for any uniqueness statement in the limited data case. A glimpse of rate results for regularized solutions completes the paper.
\end{quote}

\begin{quote}
\end{quote}

\begin{quote}
\noindent
{\small \textbf{AMS-classification (2010):}
47J06, 45Q05, 65J20
}
\end{quote}

\section{Introduction} \label{sec:intro}

Motivated by applications to spectroscopy, to the structure of solid surfaces and to nano-structures (see, e.g.,~\cite{Baumeister91,Dai13,Fukuda10,Schleicher83}) the inverse problem of \textsl{deautoconvolution}, which means that a function $x$ \textsl{with compact support} is to be reconstructed from its autoconvolution $y=x*x$, has been considered for the \textsl{one-dimensional case} extensively in the literature of the past decades. Ill-posedness, uniqueness and ambiguity as well as regularization of the deautoconvolution problem for a \textsl{real-valued function} with
compact support had been first analyzed in \cite{GorHof94}. Subsequent studies in this direction can be found in \cite{BuerFlem15,BuerHof15,BuerMat17,ChoiLan05,DaiLamm08,FleiGorHof99,FleiHof96,Janno00}. After the turn of the millennium,
the one-dimensional deautoconvolution problem for a \textsl{complex-valued function} with compact real support became of interest for modern methods of ultrashort laser pulse characterization, and we refer in this context to the article
\cite{Gerth14} as well as to the further mathematical studies in \cite{Anzengruber16,BuerFlemHof16,Flemming18}.

The object of research in this article is to present an ensemble of results for the deautoconvolution problem in the multi-dimensional case in an $L^2$-setting. We are going to extend, with respect to the reconstruction of real functions with $n$ real variables, assertions on \textsl{uniqueness, ambiguity} and \textsl{ill-posedness} that previously had been proven in the one-dimensional case. We also complement and generalize findings of our recent paper \cite{DHW22}, where such results have been stated for the two-dimensional case. Our focus is on the reconstruction of a square integrable real function $x=x(t)$ with  $t=(t_1,t_2,...,t_n)^T \in  \R^n$  of $n \ge 2$ variables with support in the unit $n$-cube $[0,1]^n$ from its autoconvolution $[x*x](s)=y(s)$ with $s=(s_1,s_2,...,s_n)^T\in \R^n$. In this context, the elements $x$ and $y$ both can be considered as tempered distributions with compact support, where $\supp(\cdot)$ is regarded as the essential support with respect to the Lebesgue measure $\lambda$ in $\R^n$. Precisely, we consider $x$ as an element of the real Hilbert space $L^2(\R^n)$ with $\supp(x) \subseteq [0,1]^n$. For short, in such a case we write $x \in L^2([0,1]^n)$ by taking into
account that $x(t)$ is assumed to be zero for $t \in \R^n \setminus [0,1]^n$.
It is well-know that, for the convolution of two functions $f$ and $g$ with $f,g \in L^2(\R^n)$ and compact supports, it holds that $f*g \in L^2(\R^n)$ as well as
\begin{equation} \label{eq:convineq}
	\supp(f*g) \subseteq \supp(f)+\supp(g)\,.
\end{equation}
Here, we use the arithmetic sum $A+B$ of two subsets $A$ and $B$ of $\R^n$ defined as $$A+B=\{a+b \in \R^n:\,a \in A,\, b \in B\}.$$
As a consequence of \eqref{eq:convineq} we have for $x \in L^2(\R^n)$ with $\supp(x) \subseteq [0,1]^n$ that $y=x*x \in L^2(\R^n)$ with $\supp(x*x) \subseteq [0,2]^n$, or in other words that $y \in L^2([0,2]^n)$.

The simplest application of our deautoconvolution problem in $n$ dimensions is the recovery of the square integrable \textsl{density function} $x$ of an $n$-dimensional random variable $\mathfrak{X}$
with support in the unit $n$-cube $[0,1]^n$
from observations of the density function $y=x*x$ of the $n$-dimensional random variable $\mathfrak{Y}:= \mathfrak{X}_1 +\mathfrak{X}_2 $, where $\mathfrak{X}, \mathfrak{X}_1$ and $\mathfrak{X}_2$ are assumed to be of i.i.d.~type. Then by definition the function $x$ obeys
the conditions $x(t)\ge 0$ a.e.~on $[0,1]^n$ and $\int_{\R^n} x(t)\,dt= \int_{[0,1]^n} x(t)\,dt= 1$.

The inverse problem of deautoconvolution is equivalent to the solution of a quadratic-type operator equation
\begin{equation} \label{eq:opeq}
	F(x)=y
\end{equation}
with the nonlinear Volterra integral operator $F: \mathcal{D}(F) \subseteq X \to Y$ mapping between the real Hilbert spaces $X:=L^2([0,1]^n)$ and $Y$ with norms $\|\cdot\|_X$ and $\|\cdot\|_Y$, respectively, and having the domain  $\mathcal{D}(F)$. Here, the nonlinear operator $F$ possesses the convolution integral form
\begin{equation} \label{eq:F}
	[F(x)](s):= [x*x](s)= \int \limits _{\R^n} x(s-t)\,x(t) dt \quad \quad (s,t \in \R^n)\,.
\end{equation}

To ease the notation in the sequel, we will make use of the following abbreviations of $n$-dimensional cuboids and cubes. If $s,t \in \R^n$ are given, we denote by
\[
\left[s,t\right]^n := \left[s_1, t_1\right]\times ... \times \left[s_n,t_n\right]
\]
the corresponding $n$-cuboids spanned by $s$ and $t$. Clearly, if $s_j >t_j$ for some $j \in \{1,...,n\}$, then $\left[s,t\right] = \emptyset$. Note that - with a slight abuse of this notation - for $s,t \in \R$ we also write $\left[s,t\right]^n$ for the $n$-cube of the form $\left[s,t\right]\times ... \times \left[s,t\right]$.

In this paper, we will distinguish two data situations. First we consider the \textsl{full data case} with $X:=L^2([0,1]^n)$, $Y:=L^2([0,2]^n)$ and forward operator $F$ as
\begin{equation} \label{eq:F1}
	[F(x)](s):=  \int\limits_{\left[\max(s-1,0), \min(s,1)\right]^n} x(s-t)\,x(t) \,\mathrm d t\,,
\end{equation}
where $y(s)=[F(x)](s)$ is observable for all $s \in [0,2]^n$, which implies due to \eqref{eq:convineq} that all relevant information about $x*x$ is available, but in practice based on noisy data $y^\delta \in Y$ with noise level $\delta>0$ and deterministic noise model
\begin{equation} \label{eq:noise}
	\|y-y^\delta\|_Y \le \delta.
\end{equation}
Secondly, we are treating again with noise model \eqref{eq:noise} the \textsl{limited data case} with \linebreak $X=Y:=L^2([0,1]^n)$ and forward operator $F$ as
\begin{equation} \label{eq:F2}
	[F(x)](s):=  \int\limits_{\left[0,s\right]^n}x(s-t)\,x(t) \,\mathrm d t\,.
\end{equation}
Here, $y(s)=[F(x)](s)$ is only available for $s$ on the unit $n$-cube $[0,1]^n$.
Since here the scope of the data is only $1/2^n$ compared to the full data case, the chances of accurately recovering $x$ from noisy observations of $y$ are decreasing more and more in the limited data case as $n$ gets larger.
In contrast to the full data case, where we assume in the sequel that $\mathcal{D}(F)=X=L^2([0,1]^n)$, we focus in the limited data case also on the domain $\mathcal{D}(F)=\mathcal{D}^+$ defined as
\begin{equation} \label{eq:domainnonnegative}
	\mathcal{D}^+:=\{x \in X=L^2([0,1]^n):\, x \ge 0 \;\; \mbox{a.e.~on}\;[0,1]^n\}\,.
\end{equation}
This set $\mathcal{D}^+$ collects the non-negative functions from $L^2([0,1]^n)$ and contains as a subset the square integrable density functions with support in the unit $n$-cube.

It is well known that, in an $L^2$-setting, the nonlinear autoconvolution operator $F$ is weakly sequentially continuous and \textsl{non-compact}, but possesses everywhere a \textsl{compact} Fr\'echet derivative of linear convolution type. This is true for
$F: X=L^2([0,1]^n) \to Y$ in both data cases $Y=L^2([0,2]^n)$ and $Y=L^2([0,1]^n)$, where the Fr\'echet derivative $F^\prime(x): X \to Y$ for all $x \in X$ attains the form
\begin{equation} \label{eq:Frechet}
	F^\prime(x)\, h =2 \,x*h  \qquad (h \in X).
\end{equation}
With this Fr\'echet derivative, the operator $F$ satisfies the nonlinearity condition
\begin{equation} \label{eq:nlcond}
	\|F(\tilde x)-F(x)-F^\prime(x)(\tilde x-x)\|_Y =\|F(\tilde x-x)\|_Y \le \|\tilde x-x\|_X^2 \qquad (\tilde x,x \in X)
\end{equation}
such that the \textsl{degree of nonlinearity} of $F$ in the sense of \cite[Def.~1]{HofSch94} is $(0,0,2)$. Note that a degree of nonlinearity $(c_1,0,0)$ with $0<c_1 \le 1$ of \textsl{tangential cone condition-type} has even not been shown for $F$ in the one-dimensional situation of $n=1$.

For any function $x \in L^2([0,1]^n)$ the autoconvolution products $F(x)=x * x$ and $F(-x)=(-x)*(-x)$ coincide for both data cases. However, it is of interest whether for $y=x*x$ the elements $x$ and $-x$ are the only solutions of equation \eqref{eq:opeq} with \linebreak $F: L^2([0,1]^n) \to L^2([0,2]^n) $ from \eqref{eq:F1} in the full data case or not. In the limited data case, for $F: \mathcal{D}^+ \subset L^2([0,1]^n) \to L^2([0,1]^n)$, it is of interest whether the solution $x$ is under non-negativity constraints even uniquely determined. Based on different versions of the Titchmarsh convolution theorem some answers to both questions are given in Section~\ref{sec:unique} below. Before that, we recall in Section~\ref{sec:prel} some basic assertions on convolution
in form of three lemmas and the definitional concept of local ill-posedness for nonlinear operator equations. In Section~\ref{sec:ill-posed}, it will be shown that the $n$-dimensional deautoconvolution problem leads in both data cases to operator equations \eqref{eq:opeq}, which are locally ill-posed everywhere. This requires the use of some kind of regularization in order to construct stable approximate solutions. Even though a detailed study on regularization of the problem is beyond the scope of this manuscript, we briefly report on some error norms and rate results for regularized solutions occurring in a numerical case study in Section \ref{sec:discuss}. There, we restrict ourselves for simplicity to the classical variant of quadratic Tikhonov regularization for nonlinear operator equations along the lines of the seminal paper \cite{EnglKunNeu89} and best possible regularization parameters. For a more detailed numerical study in the two-dimensional case we refer to  \cite{DHW22}.

\section{Preliminaries} \label{sec:prel}

Unfortunately, the formula \eqref{eq:convineq} concerning the support of the convolution function $f*g$ is an inclusion and not an equation. However, for $n=1$ and functions $f,g \in L^2(\R)$ with compact supports, which are not identically zero a.e., one can formulate an equation for the minima (smallest values) of the supports as
\begin{equation} \label{eq:infsum}
	\min \supp (f*g)\,=\, \min \supp(f)+ \min \supp(g)\,,
\end{equation}
which is a consequence of the Titchmarsh convolution theorem from \cite{Titchmarsh26}. Based on \eqref{eq:infsum} it could be shown in \cite[Theorem~1]{GorHof94} that the one-dimensional deautoconvolution problem has a uniquely determined solution in the limited data case under non-negativity constraints. By the same argument it could be shown in \cite[Theorem~4.2]{Gerth14} that $x$ and $-x$ are the only solutions in the full data case of the one-dimensional deautoconvolution problem.
An extension of those uniqueness and twofoldness results to the $n$-dimensional deautoconvolution problem require extensions of Titchmarsh's theorem to the multi-dimensional case, and we recall two versions of such extension by the following two lemmas.

The first lemma goes back to Lions (cf.~\cite{Lions51,Lions53}) and replaces $\min \supp(f)$, the \textsl{support minimum}  occurring in \eqref{eq:infsum} for $n=1$, with the \textsl{convex support} occurring in Lemma~\ref{lem:tit2} for general $n \in \N$. Here, $\csupp(f)$ denotes the convex hull of $\supp(f)$, i.e., the smallest closed convex set outside which the function $f$ vanishes a.e.~on $\R^n$.

\begin{lemma} \label{lem:tit2}
	Let the functions $f,g \in L^2(\R^n)$ with $n \in \N$ have compact supports $\supp(f)$ and $\supp(g)$. Then we have for the convolution that $f*g \in L^2(\R^n)$ and that the equation
	\begin{equation} \label{eq:sum2}
		\csupp(f*g)\,=\, \csupp(f)+\csupp(g)
	\end{equation}
	holds true. In the special case that $\supp(f*g)=\emptyset$, i.e., the function $f*g$ vanishes a.e.~on $\R^n$, then we have that at least one of the sets $\supp(f)$ or $\supp(g)$ is the empty set, which means that
	the underlying function $f$ or $g$ vanishes a.e.~on $\R^n$.
\end{lemma}

Lemma~\ref{lem:tit2} will allow us to prove the twofoldness assertion for the full data case of the multi-dimensional deautoconvolution problem in Theorem~\ref{thm:tit2} below.

\smallskip

We also present an extension of the Titchmarsh convolution theorem to the multi-dimensional case by using Mikusi\'{n}ski's $n$-\textsl{hyperpyramid} technique adapted to our situation as Lemma~\ref{lem:tit3}, and we refer in this context to
\cite[Theorem VIII]{Miku61}\,.

\begin{lemma} \label{lem:tit3}
	Let us introduce for $\gamma \ge 0$ the $n$-hyperpyramids $$\Delta(\gamma):=\{(t_1,t_2,...,\,t_n)^T\in \R^n:\;0 \le t_1,\;0 \le t_2,\;...\,,\;0 \le t_n, \; t_1+t_2+...+t_n \le \gamma\}$$ in $\R^n$\,.
	For functions $f,g \in L^2(\R^n)$ with compact supports $\supp(f)$ and $\supp(g)$ covered by $[0,\infty)^n$,  we conclude from
	$$[f*g](s)=\int \limits_{\R^n} f(s-t)\,g(t)\,dt\,=\, 0 \quad\mbox{a.e.\, for} \quad s \in \Delta(\gamma)\quad (\gamma \ge 0)$$
	that there are numbers $\gamma_1,\gamma_2 \ge 0$ with $\gamma_1+\gamma_2 \ge \gamma$ such that
	$$ f(t)\, = \, 0 \quad \mbox{a.e.\, for}\;\; t \in \Delta(\gamma_1) \quad \mbox{and} \quad g(t) \,=\, 0 \quad \mbox{a.e.\, for}\;\; t \in \Delta(\gamma_2)\,.  $$
\end{lemma}

Lemma~\ref{lem:tit3} will be used in Theorem~\ref{thm:unique} below to prove that for the limited data case of the multi-dimensional deautoconvolution problem under non-negativity constraints the solution is uniquely determined.

\smallskip

As an inverse problem the operator equation \eqref{eq:opeq} with forward operator \eqref{eq:F1} mapping from the real Hilbert space $X=L^2([0,1]^n)$ to the Hilbert space $Y=L^2([0,2]^n)$ in the full data case of multi-dimensional
deautoconvolution tends to be ill-posed. A probably stronger ill-posedenss phenomenon is to be expected for the limited data case under non-negativity constraints where $F: \mathcal{D}(F) \subset X=L^2([0,1]^n) \to Y=L^2([0,1]^n)$ and $\mathcal{D}(F)=\mathcal{D}^+$ characterize the forward operator.
For a precise theoretical verification of the ill-posedness phenomenon we adopt the \textsl{concept of local ill-posedness} for nonlinear operator equations, and we recall this concept by the following definition (cf., e.g., \cite[Def.~1.1]{HofSch98}).

\begin{definition}  \label{def:posed}
	An operator equation $F(x)=y$ with nonlinear forward operator \linebreak $F:\mathcal{D}(F) \subseteq X \to Y$  mapping between the Hilbert spaces $X$ and $Y$ with domain $\mathcal D(F)$ is called \textsl{locally ill-posed} at a solution point
	$\xdag \in \mathcal{D}(F)$ if there exist, for
	all closed balls $\overline{\mathcal{B}_r(\xdag)}$ with radius $r>0$ and center $\xdag$,
	sequences $\{x_k\} \subset \overline{\mathcal{B}_r(x^\dagger)} \cap \mathcal{D}(F)$ that satisfy the condition
	$$\|F(x_k)-F(x^\dagger)\|_Y \to 0 \, , \quad \mbox{but} \quad
	\|x_k-x^\dagger\|_X \not \to 0 \, ,
	\quad \mbox{as} \quad k \to \infty\;.$$
	Otherwise, the operator equation is called  \textsl{locally well-posed} at $\xdag$.
\end{definition}

\smallskip

For $n=1$, \textsl{local ill-posedness everywhere} on the non-negativity domain $$\mathcal{D}(F)=\{x \in X=L^2([0,1]):\, x \ge 0 \;\; \mbox{a.e.~on}\;[0,1]\}$$   was proven for the deautoconvolution problem in the
limited data case in \cite[Lemma~6]{GorHof94}. We will extend this result to the multi-dimensional situation below in Theorem~\ref{thm:ill_limited}.

Local ill-posedness everywhere on $L^2([0,1])$ could also be shown for the full data case of deautoconvolution and $n=1$ in \cite[Prop.~2.3]{FleiHof96} by perturbing the solution with an appropriate sequence of square integrable real functions, which is weakly convergent in $L^2([0,1])$. By considering such sequences as `rank one perturbations' we can also show local ill-posedness everywhere on $L^2([0,1]^n)$ in the multi-dimensional situation of deautoconvolution
with full data. For preparation we present here the following lemma, the proof of which is an immediate consequence of Lebesgue's dominated convergence theorem.

\begin{lemma} \label{lem:rank1}
	Let $\{h_k\}_{k=1}^\infty\subset L^2([0,1])$ be a sequence of real functions of one real variable, which is weakly convergent to zero, i.e.~$h_k \rightharpoonup 0$ in $L^2([0,1])$ as $k \to \infty$. Then we have for arbitrary real functions $f$ of $n$ real variables with $f \in L^2([0,1]^n)$ that the sequence  $\{f_k=f+h_k\}_{k=1}^\infty \subset L^2([0,1]^n)$
	defined as
	$$ f_k(t_1,t_2,...,t_n):=f(t_1,t_2,...,t_n)+h_k(t_1)\quad \qquad  ((t_1,t_2,...,t_n)^T \in \R^n,\;k \in \N)$$
	is weakly convergent to $f$, i.e.~$f_k \rightharpoonup f$ in $L^2([0,1]^n)$ as $k \to \infty$. In this context, we have that $\|f_k-f\|_{L^2([0,1]^n)}=\|h_k\|_{L^2(0,1)}$ for all $k \in \N$.
\end{lemma}

\section{New assertions on twofoldness and uniqueness for the multi-dimensional deautoconvolution problem}\label{sec:unique}

\subsection{Results for the full data case}

\begin{definition} \label{def:full}
	For given $y \in L^2([0,2]^n)$, we call $\xdag \in L^2([0,1]^n)$ a \textsl{solution} to the operator equation \eqref{eq:opeq} with $F: L^2([0,1]^n) \to  L^2([0,2]^n)$ according to \eqref{eq:F1} \textsl{in the full data case}
	if it satisfies the condition
	\begin{equation} \label{eq:full}
		[\xdag * \xdag](s)=y(s) \quad \mbox{a.e. for} \quad  s\in [0,2]^n\,.
	\end{equation}
\end{definition}

Lemma~\ref{lem:tit2} allows us to prove the following theorem on solution twofoldness in the full data case of multi-dimensional deautoconvolution.

\begin{theorem} \label{thm:tit2}
	If, for given $y \in L^2([0,2]^n)$, the function $\xdag \in L^2([0,1]^n)$ is a solution in the full data case  in the sense of Definition~\ref{def:full},
	then $\xdag$ and $-\xdag$ are the only solutions in this sense.
\end{theorem}
\begin{proof}
	Let $\xdag \in L^2([0,1]^n)$ be a solution in the full data case in the sense of Definition~\ref{def:full} and consider an arbitrary function $h \in L^2([0,1]^n)$ such that $\xdag+h$ is also a
	solution in the full data case in the sense of Definition~\ref{def:full}. This means that $[\xdag*\xdag](s)=[(\xdag+h)*(\xdag+h)](s)$ a.e. for $s \in [0,2]^n$, which implies that $[(\xdag+h)*(\xdag+h)-\xdag*\xdag](s)=[h*(2\xdag+h)](s)=0$ a.e. for $s \in [0,2]^n$. By setting $f:=h$ and $g:=2\xdag+h$ we can apply Lemma~\ref{lem:tit2}. Taking into account that $\supp(h*(2\xdag+h)) \subseteq[0,2]^n$, we then have $[h*(2\xdag+h)](s)=0$ a.e. for $ s \in \R^n$, or in other words $\supp(h*(2\xdag+h))=\emptyset$ and consequently $\csupp(h*(2\xdag+h))=\emptyset$. This implies, due to equation \eqref{eq:sum2}, that either $\supp(h)=\emptyset$ or $\supp(2\xdag+h)=\emptyset$ is true.
	On the one hand, $\supp(h)=\emptyset$ leads to the solution $\xdag$ itself, whereas on the other hand $\supp(2\xdag+h)=\emptyset$ leads to $[2\xdag+h](t)=0$ a.e. for $t \in [0,1]^n$
	and thus with $h=-2\xdag$ to the second solution $-\xdag$. Alternative solutions are therefore excluded. This proves the theorem.
\end{proof}

\subsection{Results for the limited data case}

\begin{definition} \label{def:limited}
	For given $y \in L^2([0,1]^n)$, we call $\xdag \in L^2([0,1]^n)$  a \textsl{solution} to the operator equation \eqref{eq:opeq} with
	$F: L^2([0,1]^n) \to L^2([0,1]^n)$ according to \eqref{eq:F2} \textsl{in the limited data case} if it satisfies the condition
	\begin{equation} \label{eq:limited}
		[\xdag * \xdag](s)=y(s) \quad \mbox{a.e. for} \quad  s\in [0,1]^n\,.
	\end{equation}
	If moreover such solution satisfies the condition  $\xdag \in \mathcal{D}^+$ with  $\mathcal{D}^+$ from \eqref{eq:domainnonnegative}, then we call it a \textsl{non-negative solution in the limited data case}.
\end{definition}

\smallskip

For solutions  $\xdag$ in the limited data case in the sense of Definition~\ref{def:limited} it is important whether the condition $0 \in \supp(\xdag)$ or its counterpart $0 \notin \supp(\xdag)$ is fulfilled.
In this context, $0 \in \supp(\xdag)$ means that for any ball $B_r(0)$ around the origin with arbitrary small radius $r>0$ there exists a set $M_r \subset B_r(0) \cap [0,1]^n$ with Lebesgue measure $\lambda(M_r)>0$ such that
$\xdag(t) \not=0$ a.e. for $t \in M_r$. Vice versa, for $0 \notin \supp(\xdag)$ we have some sufficiently small radius $r>0$ such that $\xdag(t)=0$ a.e. for $t \in B_r(0) \cap [0,1]^n$.

\smallskip

In a first step we generalize by Theorem~\ref{thm:nonunique} those aspects that had been fixed for $n=1$ in \cite[Theorem~1]{GorHof94} concerning the strong non-injectivity of the autoconvolution operator in the limited data case
to the multi-dimensional situation with arbitrary $n \in \N$.

\begin{theorem} \label{thm:nonunique}
	If, for given $y \in L^2([0,1]^n)$, the function $\xdag \in L^2([0,1]^n)$ is a solution in the limited data case in the sense of Definition~\ref{def:limited} that fulfills the condition
	\begin{equation} \label{eq:notzero1}
		0 \notin \supp(\xdag)\,,
	\end{equation}
	then there exist infinitely many other solutions ${\hat x}^\dagger \in L^2([0,1]^n)$ in this sense.
\end{theorem}
\begin{proof}
	Under the condition \eqref{eq:notzero1} there is some $0<\varepsilon<1/2$ such that $\xdag(t)=0$ a.e. for $t \in [0,\varepsilon]^n$.
	Now there exist infinitely many $h \in L^2([0,1]^n)$ such that
	\begin{equation} \label{eq:hnull}
		h(t)\,=\,0 \quad \mbox{a.e. for} \;\; t \in [0,1]^n \setminus [1-\varepsilon,1]^n\,.
	\end{equation}
	For any such $h$ we have
	\[
	[(2\xdag+h)*h](s) = \int_{[0,s]^n} (2\xdag + h)(s-t) \,h(t) \,\mathrm d t = \int_{[1-\varepsilon,s]^n} (2\xdag + h)(s-t)\, h(t) \,\mathrm d t.
	\]
	But for $s \in [0,1]^n$ and $t \in [1-\varepsilon,s]^n$, we have component-wise that
	\[
	0 \leq s_i - t_i \leq 1-t_i \leq \varepsilon < 1-\varepsilon
	\]
	due to $\varepsilon < \frac12$, so that $(2\xdag + h)(s-\cdot) = 0$ a.e. for $[1-\varepsilon,s]^n$. Therefore
	\[
	[(2\xdag+h)*h](s) = 0  \quad \mbox{a.e. for} \quad  s \in [0,1]^n,
	\]
	which implies
	$$[\left(\xdag + h\right) * \left(\xdag + h\right)] (s) = \xdag * \xdag(s) + [(2\xdag+h)*h](s)=y(s) \quad \mbox{a.e. for} \quad  s \in [0,1]^n\,.$$
	This yields the claim.\end{proof}

Now we are ready to formulate and to prove with the following theorem a main new result of this paper, which extends the solution uniqueness assertion for the limited data case under non-negativity constraints published for $n=1$
in \cite[Theorem~1]{GorHof94} to the multi-dimensional situation with arbitrary $n \in \N$.
The proof of this theorem is based on Mikusi\'{n}ski's $n$-hyperpyramid technique introduced above by Lemma~\ref{lem:tit3}.

\begin{theorem} \label{thm:unique}
	If, for given $y \in L^2([0,1]^n)$, the function $\xdag \in L^2([0,1]^n)$ is a non-negative solution in the limited data case in the sense of Definition~\ref{def:limited} that fulfills the condition
	\begin{equation} \label{eq:notzero3}
		0 \in \supp(\xdag)\,,
	\end{equation}
	then $\xdag$ is the uniquely determined non-negative solution in this case.
\end{theorem}
\begin{proof}
	First we will show that under the condition \eqref{eq:notzero3} the non-negative solution $\xdag(t)$ is uniquely determined a.e.~for $t \in \Delta(1)$. Namely, supposed that there exists a function $h \in L^2([0,1]^n)$ with $\xdag+h \ge 0$  satisfying the equation
	\begin{equation} \label{eq:gapclose}
		[(\xdag+h) * (\xdag+h)](s)=y(s) \quad \mbox{a.e. for} \quad  s \in [0,1]^n\,,
	\end{equation}
	we would have that
	\begin{equation} \label{eq:hgapclose}
		[h*(2\xdag+h)](s)=0\quad \mbox{a.e. for} \quad  s \in [0,1]^n.
	\end{equation}
	Because of $[0,1]^n \supset \Delta(1)$, Lemma~\ref{lem:tit3} applies with $f:=h, \; g:=2\xdag+h$ and $\gamma=1$. Obviously, we have  $\gamma_2=0$ due to the fact that $[2\xdag+h](t) \ge \xdag(t)$ a.e. for $t \in [0,1]^n$,
	which implies together with condition \eqref{eq:notzero3} that $0 \in \supp(2\xdag+h)$. Then we find as a consequence of $\gamma_1+\gamma_2 \ge \gamma$ that $\gamma_1 \ge 1$ must hold, which yields
	$h(t)=0$ a.e. for $t \in \Delta(1)$.
	
	In a second step of the proof we show that also perturbations $h \in L^2([0,1]^n)$ with $\xdag+h \ge 0$ and  $\supp(h) \subseteq \overline{[0,1]^n\setminus \Delta(1)}$  are
	only possible if $h$ is the zero function almost everywhere on $[0,1]^n \cap \Delta(2)$.
	Now assume, for such function $h$, that it obeys the condition \eqref{eq:gapclose} and consequently \eqref{eq:hgapclose}.
	From \eqref{eq:hgapclose} we derive that
	\begin{equation} \label{eq:igapclose}
		[(-2\xdag)*h](s)= [h*h](s) \quad \mbox{a.e. for} \quad  s \in [0,1]^n\,,
	\end{equation}
	which allows us to apply Lemma~\ref{lem:tit3} with $f:=-2\xdag$, $g:=h$, $f*g=h*h$ and the associated values $\gamma_1, \gamma_2$ and $\gamma$, respectively.
	Evidently, we have $$\supp(h*h) \subseteq 2 \supp(h) \subseteq \overline{[0,2]^n\setminus \Delta(2)}$$ and thus $\gamma=2$. This yields $\gamma_2=2$ and hence $h=0$ a.e.~on $[0,1]^n \cap \Delta(2)$, because $\gamma_1=0$
	as a consequence of condition \eqref{eq:notzero3}. Now, for $n=2$ the proof is complete, because of $[0,1]^2 \subset [0,1]^2 \cap \Delta(2)$. For $n>2$, however we must repeat the second step in an analog manner $m$ times until $2^m \ge n$
	such that $h=0$ a.e.~on $[0,1]^n \cap \Delta(2^m) \supseteq [0,1]^n \cap \Delta(n) = [0,1]^n$. Then the proof is complete.
\end{proof}

\section{Ill-posedness phenomena} \label{sec:ill-posed}

For nonlinear inverse problems modelled by operator equations \eqref{eq:opeq} in Hilbert spaces, the character and strength of ill-posedness may be a local property and may depend on nonlinearity conditions
of the forward operator $F$, see for discussions and examples of the articles \cite{Hof94,HofPla18,HofSch94}. Therefore, the concept of local ill-posedness at a solution point $\xdag$ (see Definition~\ref{def:posed} above) applies for \eqref{eq:opeq} with the autoconvolution operator $F$ from \eqref{eq:F}. It could be proven for the one-dimensional situation that the deautoconvolution problem is \textsl{locally ill-posed everywhere} on $\mathcal{D}(F)=L^2([0,1])$
for the full data case (cf.~\cite[Prop.~2.3]{FleiHof96}) and on $\mathcal{D}(F)=\mathcal{D}^+ \subset L^2([0,1])$ with $\mathcal{D}^+$ from \eqref{eq:domainnonnegative} with $n=1$ for the limited data case (cf.~\cite[Lemma~6]{GorHof94}). The following two theorems extend the results to the multi-dimensional situation for arbitrary $n \in \N$.

\begin{theorem}\label{thm:ill_limited}
	For the limited data case of deautoconvolution, the operator equation \eqref{eq:opeq} with $X=Y=L^2([0,1]^n)$ and forward operator $F: \mathcal{D}^+ \subset X \to Y$ from \eqref{eq:F2} with non-negativity domain $\mathcal{D}^+$ from \eqref{eq:domainnonnegative} is \textsl{locally ill-posed everywhere} on $\mathcal{D}(F)=\mathcal{D}^+$.
\end{theorem}
\begin{proof}
	Let $\xdag \in \mathcal{D}^+$ be a non-negative solution in the limited data case in the sense of Definition~\ref{def:limited}. To show local ill-posedness at $\xdag$ we introduce for fixed $r>0$ the sequence $\{h_k\}_{k=3}^\infty \subset L^2([0,1]^n)$ of perturbations of the form
	$$h_k(t):=\left\{\begin{array}{ccl} k^{n/2} \, r & \mbox{for} & t \in [1-\frac{1}{k},1]^n\\ & & \\0 & \mbox{for} & t \in [0,1]^n \setminus [1-\frac{1}{k},1]^n   \end{array} \right. $$
	with $x_k:=\xdag+h_k \in \mathcal{D}^+$, $\|h_k\|_{L^2([0,1]^n)}=r$ and consequently $x_k \in \overline{{\mathcal B}_r(x^\dagger)} \cap \mathcal{D}^+$ for all $k \ge 3$.
	To complete the proof of the theorem we still need to show that
	the norm \linebreak $\|F(x_k)-F(\xdag)\|_{L^2([0,1]^n)}$ tends for all $r>0$ to zero as $k$ tends to infinity. Owing to \linebreak $F(x_k)-F(\xdag)=2\xdag*h_k+h_k*h_k$ and $\|h_k*h_k\|_{L^2([0,1]^n)}=0\,$, this rewrites as
	$$\|\xdag*h_k\|_{L^2([0,1]^n)}  \to 0 \quad \mbox{as} \quad k \to \infty\,. $$
	Evidently, for $s=(s_1,s_2,...,s_n)^T, t=(t_1,t_2,...,t_n)^T \in \R^n$, the non-negative values
	$$[\xdag*h_k](s)=\int\limits_{[0,s]^n} h_k(s-t)\,\xdag(t)\,\mathrm dt $$
	can be different from zero only for $s \in [1-\frac{1}{k},1]^n$. Using the Cauchy-Schwarz inequality and taking into account that $\xdag \in \mathcal{D}^+$
	we have for those $s \in [1-\frac{1}{k},1]^n$ the estimate
	$$ [\xdag*h_k](s)=  k^{n/2} \, r \int\limits_{\left[0,s-\left(1-\frac1k\right)\right]^n} \xdag(t)\, \mathrm d t$$
	$$ \le r\,\|\xdag\|_{L^2([0,1]^n)}.   $$
	This, however, yields
	$$\|\xdag*h_k\|_{L^2([0,1]^n)} \le r\,\,\|\xdag\|_{L^2([0,1]^n)} \left(~\int \limits_{\left[ 1-\frac{1}{k},1\right]^n} 1 \,\mathrm d s\right)^{1/2}
	=\frac{ r\,\,\|\xdag\|_{L^2([0,1]^n)}}{k^{n/2}}$$
	tending for all $r>0$ to zero as $k$ tends to infinity. This completes the proof of the theorem.
\end{proof}

\begin{theorem}\label{thm:ill_full}
	For the full data case of deautoconvolution, the operator equation \eqref{eq:opeq} with $X= L^2([0,1]^n)$, $Y=L^2([0,2]^n)$ and forward operator $F:X \to Y$ from \eqref{eq:F1} is \textsl{locally ill-posed everywhere} on $\mathcal{D}(F)=X$.
\end{theorem}
\begin{proof}
	Let $\xdag \in L^2([0,1]^n)$ be a solution in the full data case in the sense of Definition~\ref{def:full}. For showing local ill-posedness at $\xdag$ we fix $r>0$ arbitrary and introduce the sequence $\{h_k\}_{k=1}^\infty \subset L^2([0,1])$ of functions of one real variable of the form
	\begin{equation}\label{eq:h_k}
		h_k(t):=\sqrt{2}\,r\,\sin(k^2t^2) \qquad (t \in [0,1],\;k \in \N).
	\end{equation}
	For finding properties of $h_k$ and $F(h_k)=h_k*h_k$ one needs to use the Fresnel integrals
	$$S(s):=\int_0^s \sin(t^2)dt \quad \mbox{and} \quad C(s):=\int_0^s \cos(t^2)dt.$$
	For $s \in [0,\infty)$ the range of both continuous functions is covered by the interval $[0,1]$.
	One easily finds that $0.5\,r<\|h_k\|_{L^2([0,1])} < r= \lim_{k \to \infty}  \|h_k\|_{L^2([0,1])}$ and that the weak convergence
	$h_k \rightharpoonup 0$ in $L^2([0,1])$ as $k \to \infty$ takes place. The latter is a consequence of the fact that, for all $0 \le s \le 1$,
	$$ 0 \le \int_0^s h_k(t)dt= \frac{\sqrt{\pi}\,r\,S(k\sqrt{2/\pi)}s)}{k} \le \frac{\sqrt{\pi}\,r}{k}  \to 0 \quad \mbox{as} \quad k \to \infty. $$
	Now we consider the perturbed functions $x_k:=\xdag+h_k \in L^2([0,1]^n)$ defined as
	$$ x_k(t_1,t_2,...,t_n):=\xdag(t_1,t_2,...,t_n)+h_k(t_1)\quad \qquad  ((t_1,t_2,...,t_n)^T \in \R^n,\;k \in \N),$$
	with $x_k \in \overline{\mathcal{B}_r(\xdag)}$ and $\|x_k-\xdag\|_{L^2([0,1]^n)}=\|h_k\|_{L^2([0,1])} \not\to 0$ as $k \to \infty$. To complete the proof, we still have to show that
	$$\|F(x_k)-F(x^\dagger)\|_{L^2([0,2]^n)} \to 0  \quad \mbox{as} \quad k \to \infty\;.$$
	Since $F(x_k)-F(\xdag)=F^\prime( \xdag)(x_k-\xdag)+F(x_k-\xdag)$ and $x_k-\xdag \rightharpoonup 0$  as $k \to \infty$, we have
	$\lim_{k \to \infty}\|F(x_k)-F(\xdag)\|_{L^2([0,2]^n)} \le \lim_{k \to \infty}\|F(x_k-\xdag)\|_{L^2([0,2]^n)} \le \lim_{k \to \infty}\|h_k*h_k\|_{L^2([0,2])}$ by taking into account Lemma~\ref{lem:rank1} and that $F^\prime(\xdag)$ is a compact operator. To complete the proof we finally show that  $\lim_{k \to \infty}\|h_k*h_k\|_{L^2([0,2])}=0$. Owing to the properties of Fresnel integrals mentioned above, this is a consequence of $|[h_k*h_k](\xi)|\le \frac{\bar C}{k}$ for $\xi \in [0,2]$ with a uniform constant $\bar C>0$, which follows from the two formulas
	$${\textstyle \int\limits_0^s \sin(k^2(s-t)^2)\sin(k^2t^2)dt =\frac{\sqrt{\pi} k s\left( S(ks/\sqrt{\pi})\sin(k^2s^2/2)-C(ks/\sqrt{\pi})\cos(k^2s^2/2) \right)+\sin(k^2s^2)}{2k^2s}}$$
	valid for $0 \le s \le 1$, and
	$${\textstyle \int\limits_{s-1}^1 \sin(k^2(s-t)^2)\sin(k^2t^2)dt=\frac{\sqrt{\pi} k s\left( S(k(2-s)/\sqrt{\pi})\sin(k^2s^2/2)-C(k(2-s)/\sqrt{\pi})\cos(k^2s^2/2) \right)+\sin(k^2s(2-s))}{2k^2s}}$$
	valid for $1 < s \le 2$.
\end{proof}

\begin{figure}
	\begin{center}
		\includegraphics[width=7cm]{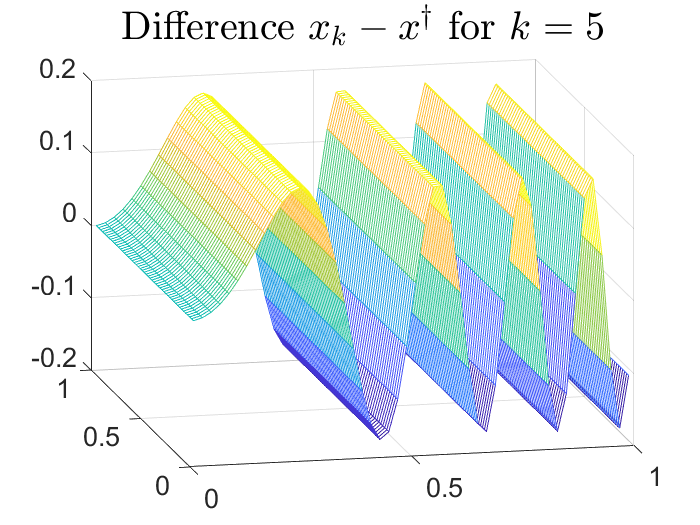}
		\includegraphics[width=7cm]{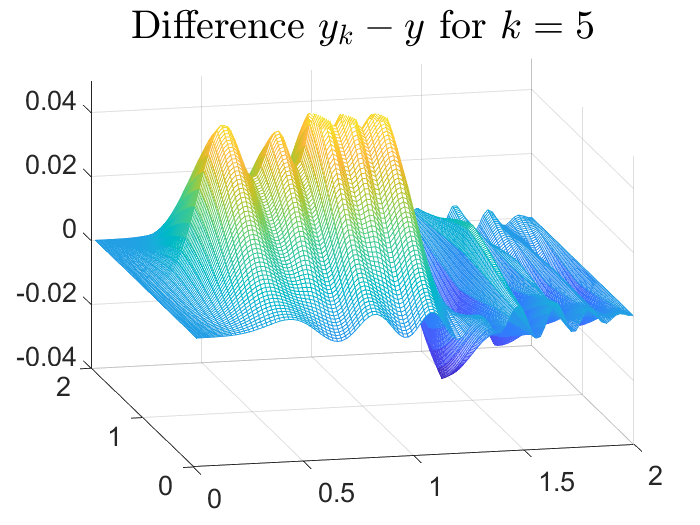}
		\includegraphics[width=7cm]{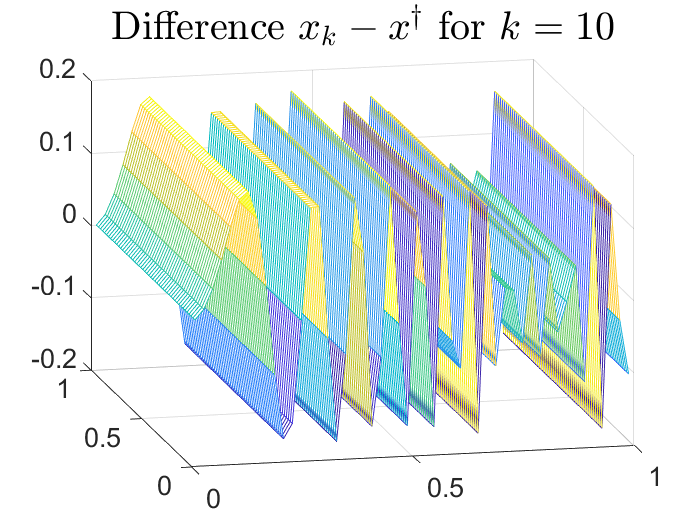}
		\includegraphics[width=7cm]{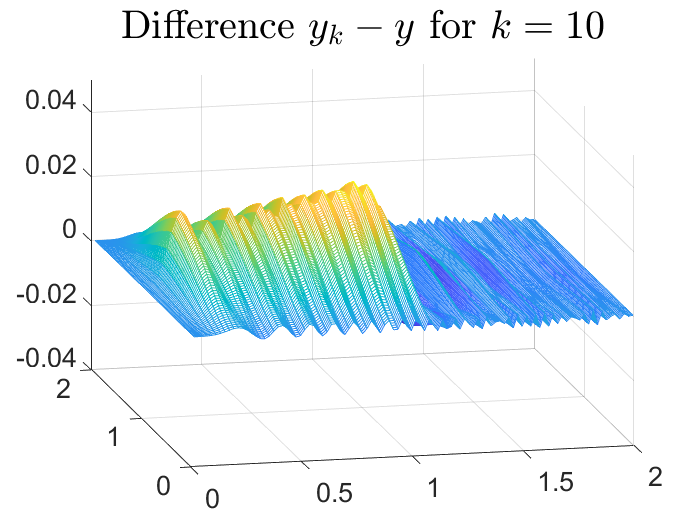}
		\includegraphics[width=7cm]{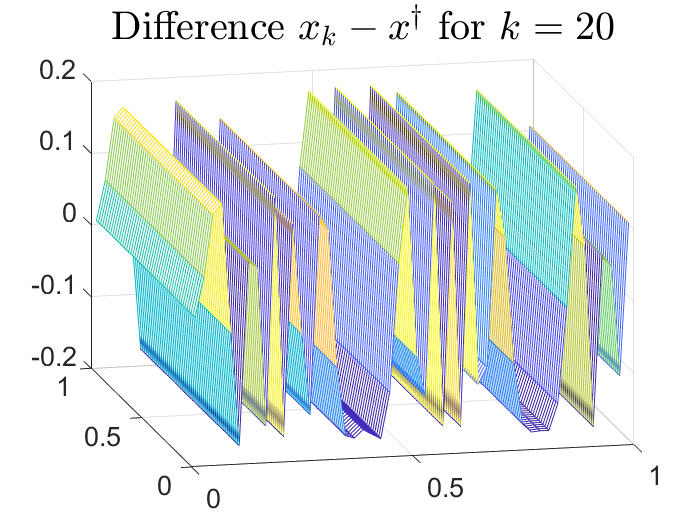}
		\includegraphics[width=7cm]{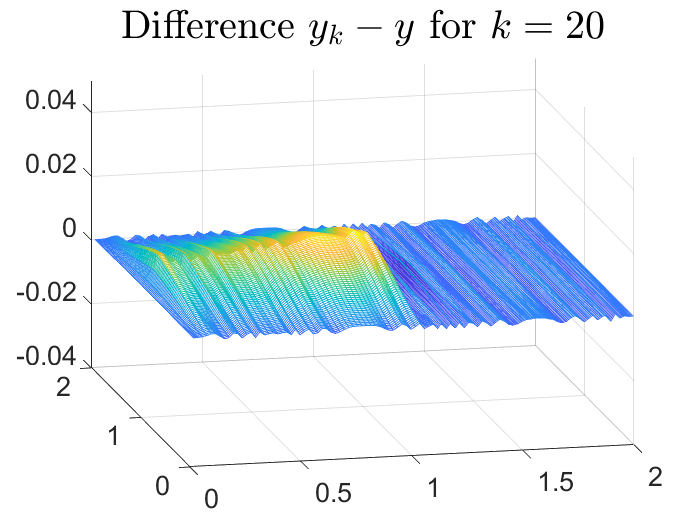}
		\caption{Development of differences $x_k-x^\dagger$ and $y_k-y$ for increasing $k$}
		\label{fig:illposed}
	\end{center}
\end{figure}

We are going to illustrate with Figure~\ref{fig:illposed} the ill-posedness phenomenon for the full data case of deautoconvolution along the lines of the ideas of the proof of Theorem~\ref{thm:ill_full}.
For this purpose we exploit as an example solution the function
\[x^\dagger(t_1,t_2)=\left[\frac{2}{3}(t_1+1)\right]\cdot\left[\frac{\pi}{2+\pi}\left(\cos\left(\left(t_2-\frac{1}{2}\right)\pi\right)+1\right)\right]\,,\]
which characterizes a factorable probability density function of a two-dimensional random vector with two uncorrelated one-dimensional components.
For the sequence introduced in \eqref{eq:h_k} we use the function $h_k(t)=\frac{\sqrt{2}}{8}\sin(k^2t^2)$, which leads to the perturbed solution $x_k(t_1,t_2)=x^\dagger(t_1,t_2)+h_k(t_1)$ that converges weakly in $L^2([0,1]^2)$ to $\xdag$ as $k \to \infty$, but not in norm as
the pictures of $x_k-\xdag$ on the left in Figure~\ref{fig:illposed} for $k=5,10$ and $20$ clearly show. However, the pictures on the right indicate convincingly the norm convergence of $y_k=x_k*x_k$ to $y=\xdag*\xdag$ in the space $L^2([0,2]^2)$.

\section{A glimpse of rate results for regularized solutions} \label{sec:discuss}

The goal of this concluding section is to mention some unexpected behavior of regularized solutions occurring in a brief case study on deautoconvolution, where in a setting analogous to \cite{EnglKunNeu89}
and \cite[Sect.~10.2]{EHN96} the regularized solutions
\begin{equation} \label{eq:Tik}
	x_\alpha^\delta \in \argmin_{x \in \mathcal D(F)} \left[ \Vert F(x) - y^\delta\Vert^2_Y + \alpha \Vert x-\bar{x}\Vert^2_X \right]
\end{equation}
are minimizers of the Tikhonov functional. For both operators \eqref{eq:F1} and \eqref{eq:F2} under consideration, the element $y^\delta \in Y$ denotes the available data satisfying \eqref{eq:noise}, $\bar{x}\in X$ is a
reference element (initial guess), and $\alpha > 0$ is a regularization parameter. Our study is reduced to the case that best possible regularization parameters $\alpha=\alpha_{opt}$ in the sense of
\begin{equation}\label{eq:alpha_opt}
	\alpha_{opt}(\delta)=\min_{\alpha>0}\|x_\alpha^\delta-x^\dagger\|_X
\end{equation}
are evaluated. From the three density functions of one real variable with supports in $[0,1]$,
\begin{equation*}
	x_1(t_1)=\frac{2(t_1+1)}{3}, \;\; x_2(t_2)=\frac{\pi}{2+\pi}\left(\cos((t_2-\frac{1}{2})\pi)+1\right),\;\;
	x_3(t_3)=\left\{\begin{array}{ll}
		\frac{5}{4}\quad\;  0\leq t_1< \tfrac12\\
		t_1\quad\; \tfrac12\leq t_1\leq 1\\
	\end{array}\right.,
\end{equation*}
we assemble two solutions $\xdag$ for the two- and three-dimensional situation of deautoconvolution as
$$\xdag(t_1,t_2)=x_1(t_1)x_2(t_2)\qquad \mbox{for}\;\;n=2 $$
and
$$\xdag(t_1,t_2,t_3)=x_1(t_1)x_2(t_2)x_3(t_3)\qquad \mbox{for} \;\;n=3, $$
which are density functions with supports $[0,1]^n$.
To the discretization level with a uniform meshwidth of $\frac{1}{50}$ in each direction, the regularized solutions $x_{\alpha_{opt}}^\delta$ according to \eqref{eq:Tik} have been calculated with a constant initial guess  $\bar{x}\equiv 0.5$ in the discretized form for $n=2,3$ and randomly generated noisy data $y^\delta \in L^2([0,2]^n)$ (full data case) as well as for $y^\delta \in L^2([0,1]^n)$ (limited data case). The relative empirical errors in \% measured in the discrete $L^2$-norm for different $\delta$, each simulated from 10 independent runs, are listed in Table~\ref{tab:1}.
The bottom line of the table contains the H\"older exponent $0<\kappa<1$ of the convergence rate $\|x_{\alpha_{opt}}^\delta-\xdag\|=\mathcal{O}(\delta^\kappa)$ as $\delta \to 0$ for the different situations, which had been
estimated by regression from the selection of $\delta$-values under consideration in the table.

\begin{table}[H]
	\caption{Relative error norms of regularized solutions }
	\label{tab:1}
	\begin{center}
		\begin{tabular}{|c|c|c|c|c|}
			\hline
			\scriptsize{relative input errors} & \multicolumn{4}{|c|}{\scriptsize{relative output errors of $x_{\alpha_{opt}}^\delta$}}\\
			$\frac{\|y^{\delta}-y^\dagger\|_Y}{\|y^\dagger\|_Y}$ & \multicolumn{4}{|c|}{$\frac{\|x_{\alpha_{opt}}^\delta-x^\dagger\|_X}{\|x^\dagger\|_X}$}\\
			\hline
			&\multicolumn{2}{|c|}{full data case}&\multicolumn{2}{|c|}{limited data case} \\
			&$n=2$ & $n=3$ & $n=2$ & $n=3$\\
			\hline
			$10\%$ & $9.85\%$ &$13.48\%$ &$17.54\%$ & $23.59\%$\\
			$8\%$ & $8.70\%$& $12.12\%$&$17.21\%$&$22.59\%$\\
			$5\%$ & $6.38\%$& $9.82\%$&$15.17\%$&$19.99\%$\\
			$2\%$ & $3.61\%$& $6.26\%$&$9.74\%$&$14.54\%$\\
			$1\%$ & $2.31\%$& $4.12\%$&$7.95\%$&$11.58\%$\\
			$0.8\%$ & $1.98\%$& $3.57\%$&$7.39\%$&$10.50\%$\\
			$0.5\%$ & $1.44\%$& $2.61\%$&$5.94\%$&$9.24\%$\\
			$0.2\%$ & $0.78\%$& $1.42\%$&$4.10\%$&$6.85\%$\\
			$0.1\%$ & $0.48\%$& $0.87\%$&$2.70\%$&$5.47\%$\\
			$0.05\%$ & $0.30\%$& $0.53\%$&$1.76\%$&$4.31\%$\\
			\hline
			\hline
			estim.~H\"older exponent $\kappa$ & $0.66$ & $0.61$ & $0.43$ & $0.32$\\
			\hline
		\end{tabular}
	\end{center}
\end{table}

An inspection of Table~\ref{tab:1} shows for both dimensions $n=2$ and $n=3$ a substantial reduction of the regularization error norms in the full data case compared to the limited data case. This is intuitively explained by the lack of data in $[0,2]^n \setminus [0,1]^n$, but even though this lack is considerably larger in dimension $n= 3$ (factor $8$) compared to $n = 2$ (factor $4$), the error norms do not fully reflect this behavior.

Based on ten different noise levels $\delta$, a rough estimation of convergence rates of the corresponding error norms as $\delta$ tends zero indicates H\"older exponents $\kappa>0.5$ in the full data case and $\kappa<0.5$ in the limited data case. However, both results cannot fully be explained by available theory. It is known from \cite{EnglKunNeu89} and \cite[Theorem~10.4]{EHN96} that a $\kappa = 0.5$ rate (i.e. $\mathcal O (\sqrt{\delta})$) is obtained under a range-type source condition $\xdag-\bar{x}=(F^\prime(\xdag))^*w$ in combination with a smallness condition on $\left\Vert w\right\Vert_Y$. On the other hand, it has been shown in \cite[Prop.~2.6]{BuerHof15} that such theory is hard to apply for the autoconvolution operator $F$ even in the one-dimensional case. To obtain rates with $\kappa > 0.5$, it is i.e. known from \cite{Neubau89} and \cite[Theorem~10.7]{EHN96}, that a rate $\mathcal O (\delta^{\frac23})$ can be obtained under the higher-order range condition $\xdag-\bar{x}=(F^\prime(\xdag))^*F^\prime(\xdag)v$ in combination with a smallness condition on $\left\Vert v\right\Vert_X$. But in view of \cite[Prop.~2.6]{BuerHof15} it is also questionable whether such a result can be applied for the autoconvolution operator $F$ at hand. In both situations, one reason seems to be fact that the compact Fr\'echet derivatives $F^\prime(x)$ carry too little information about the non-compact operator $F$. Also \textsl{variational source conditions} introduced in \cite{HKPS07} and, for example, further analyzed in \cite{HohWer22,WerHoh12} could not be successfully exploited for obtaining convergence rates in deautoconvolution. Solely in \cite[Prop.~5.1 and Cor.~5.2]{BuerFlemHof16} a convergence rate could be derived by means of variational source conditions, but only under strong sparsity assumptions on the solution $\xdag$. Nevertheless, the numerical experiment in the context of Table~\ref{tab:1} indicates the practical occurrence of H\"older convergence rates for regularized solutions to the multi-dimensional deautoconvolution problem.

\bigskip

\section*{Acknowledgment} Bernd Hofmann and Yu Deng are supported by the German Science Foundation (DFG) under the grant~HO~1454/13-1 (Project No.~453804957).

\end{document}